\documentclass{crm-article}

\usepackage{arxiv}
\usepackage{algorithm}
\usepackage{algorithmic}

\usepackage{mathtools, amssymb}
\usepackage{amsmath}

\usepackage{xcolor}
\usepackage{mathtools}

\usepackage[useregional]{datetime2}

\renewcommand{\le}{\leq}
\renewcommand{\ge}{\geq}

\renewcommand{\limits}{}
\newcommand{\argmin}{\text{arg\,min}}
\newcommand{\diam}{\text{\,diam}}

\newcommand{\cbraces}[1]{\left( #1 \right)}
\newcommand{\sbraces}[1]{\left[ #1 \right]}

\usepackage[colorinlistoftodos,bordercolor=orange,backgroundcolor=orange!20,linecolor=orange,textsize=scriptsize]{todonotes}

\graphicspath{ {pics} }

\def\<#1,#2>{\langle #1,#2\rangle}

\begin{document}





\journalVol{10}
\journalNo{1} 
\setcounter{page}{1}

\journalSection{Математические основы и численные методы моделирования}
\journalSectionEn{Mathematical modeling and numerical simulation}

\journalReceived{01.06.2016.}
\journalAccepted{01.06.2016.}

\UDC{519.8}
\title{Модификации алгоритма Frank-Wolfe в задаче поиска равновесного распределения транспортных потоков}
\titleeng{Modifications of the Frank-Wolfe algorithm in the problem of finding the equilibrium distribution of traffic flows}
\thanks{Работа Игнашина И.Н. была выполнена при поддержке ежегодного дохода ФЦК МФТИ ({целевого капитала № 5} на развитие направлений искусственного интеллекта и машинного обучения в МФТИ, \url{https://fund.mipt.ru/capitals/ck5/}).\\
Работа Ярмошика Д.В. была выполнена при поддержке Министерства науки и высшего образования Российской Федерации (госзадание), номер проекта FSMG-2024-0011    
}
\thankseng{The work of Ignashin I. is supported by the annual income of the Endowment Fund of Moscow Institute of Physics and Technology (target capital no. 5 for the development of artificial intelligence and machine learning, \url{https://fund.mipt.ru/capitals/ck5/}).\\
The work of Yarmoshik D. is supported by the Ministry of Science and Higher Education of the Russian Federation (Goszadaniye), project No. FSMG-2024-0011
}

\author[1]{\firstname{И.\,Н.}~\surname{Игнашин}}
\authorfull{Игорь Николаевич Игнашин}
\authoreng{\firstname{I.\,N.}~\surname{Ignashin}}
\authorfulleng{Igor N. Ignashin}
\email{ignashin.in@phystech.edu}
\affiliation[1]{Национальный исследовательский университет «Московский физико-технический институт»,\protect\\ Россия, 141701, г. Долгопрудный, Институтский пер., д. 9}
\affiliationeng[1]{Moscow Institute of Physics and Technology, 9 Instituskiy ave., Dolgoprudny, 141701, Russia}

\affiliation[2]{Институт проблем передачи информации РАН, \protect\\Россия, 127051, г.
Москва, Большой Каретный пер., д. 19, стр. 1}
\affiliationeng[2]{Kharkevich Institute for Information Transmission Problems, 19 Bolshoy Karenty per., Moscow, 127051, Russia}



\author[1,2]{\firstname{Д.\,В.}~\surname{Ярмошик}}
\authorfull{Демьян Валерьевич Ярмошик}
\authoreng{\firstname{D.}~\surname{Yarmoshik}}
\authorfulleng{Demyan Yarmoshik}
\email{yarmoshik.dv@phystech.edu}

\begin{abstract}

В работе приведены различные модификации алгоритма Frank-Wolfe для задачи поиска равновесного распределения потоков. В качестве модели для экспериментов используется модель Бекмана. В этой статье в первую очередь уделяется внимание выбору направления базового шага алгоритма Frank-Wolfe (FW). Будут представлены алгоритмы: Conjugate Frank-Wolfe (CFW) , Bi-conjugate Frank-Wolfe (BFW), Fukushima Frank-Wolfe (FFW). Каждой модификации соответствуют различные подходы к выбору этого направления. Некоторые из этих модификаций описаны в предыдущих работах авторов. В данной статье будут предложены алгоритмы N-conjugate Frank-Wolfe (NFW) и Weighted Fukushima Frank-Wolfe (WFFW). Эти алгоритмы являются некоторым идейным продолжением алгоритмов BFW и FFW. Таким образом, если первый алгоритм использовал на каждой итерации два последних направления предыдущих итераций для выбора следующего направления, сопряженного к ним, то предложенный алгоритм NFW использует $N$ предыдущих направлений. В случае же Fukushima Frank-Wolfe~--- в качестве следующего направления берется среднее от нескольких предыдущих направлений. Соответственно этому алгоритму предложена модификация WFFW, использующая экспоненциальное сглаживание по предыдущим направлениям. Для сравнительного анализа были проведены эксперименты с различными модификациями на нескольких наборах данных, представляющих городские структуры и взятых из общедоступных источников. За метрику качества была взята величина относительного зазора. Результаты экспериментов показали преимущество алгоритмов, использующих предыдущие направления для выбора шага, перед классическим алгоритмом Frank-Wolfe.  Кроме того, было выявлено улучшение эффективности при использовании болee двух сопряженных направлений. Например, на многих датасетах модификация 3-conjugate FW сходилась наилучшим образом. Кроме того, предложенная модификация WFFW зачастую обгоняла FFW и CFW, хотя и проигрывала модификациям NFW.

\end{abstract}




\begin{abstracteng}

The paper presents various modifications of the Frank-Wolfe algorithm in the equilibrium traffic assignment problem. The Beckman model is used as a model for experiments. In this article, first of all, attention is paid to the choice of the direction of the basic step of the Frank-Wolfe algorithm. Algorithms will be presented: Conjugate Frank-Wolfe (CFW), Bi -conjugate Frank-Wolfe (BFW), Fukushima Frank-Wolfe (FFW). Each modification corresponds to different approaches to the choice of this direction. Some of these modifications are described in previous works of the authors. In this article, following algorithms will be proposed: N-conjugate Frank-Wolfe (NFW), Weighted Fukushima Frank-Wolfe (WFFW). These algorithms are some ideological continuation of the BFW and FFW algorithms. Thus, if the first algorithm used at each iteration the last two directions of the previous iterations to select the next direction conjugate to them, then the proposed algorithm NFW is using more than $N$ previous directions. In the case of Fukushima Frank-Wolfe, the average of several previous directions is taken as the next direction. According to this algorithm, a modification WFFW is proposed, which uses a exponential smoothing from previous directions. For comparative analysis, experiments with various modifications were carried out on several data sets representing urban structures and taken from publicly available sources. The relative gap value was taken as the quality metric. The experimental results showed the advantage of algorithms using the previous directions for step selection over the classic Frank-Wolfe algorithm.  In addition, an improvement in efficiency was revealed when using more than two conjugate directions. For example, on various datasets, the modification 3FW showed the best convergence. In addition, the proposed modification WFFW often overtook FFW and CFW, although performed worse than NFW.

\keywordeng{Conjugate Frank-Wolfe}
\keywordeng{Weighted Fukushima Frank-Wolfe}
\keywordeng{N-conjugate Frank-Wolfe}
\end{abstracteng}

\maketitle

\paragraph{Введение}

Задача поиска равновесного распределения потоков является одним из фундаментальных аспектов транспортного планирования. На данный момент, методологии, основанные на машинном обучении, не достигли эффективности, сравнимой с подходами на основе теоретико-игровых моделей \cite{datadriven2022}. Поэтому исследования продолжаются в рамках модели Бекмана. 
Проводятся исследования различных алгоритмов для решения задачи поиска равновесного распределения потоков. Так, в работе \cite{babazadeh2020} представлен алгоритм, использующий пути в явном виде, потому требующий больше памяти, но достаточно эффективно работающий по времени. Также эффективными в данной задаче являются алгоритмы, основанные на использовании потоков по дорогам, к которым относится алгоритм Frank-Wolfe.

Алгоритм Frank-Wolfe нашел применение во многих областях, в том числе в транспортной логистике и телекоммуникациях. Этот метод оптимизации решает задачи математического программирования с ограничениями, но, несмотря на свою эффективность, его основной недостаток заключается в быстром замедлении сходимости после первых десятков итераций, что создает необходимость искать модификации этого алгоритма.


Алгоритм Frank-Wolfe (FW) был предложен в работе \cite{Frank1956}. В статье \cite{pedregosa2020linearly} проведен современный анализ сходимости алгоритма FW и получены оценки, которые могут быть использованы в контексте задачи, рассматриваемой в настоящей статье. 

 
Одна из наиболее эффективных модификаций алгоритма FW~--- это Conjugate Frank-Wolfe (CFW). Данная модификация использует для выбора направления шага предыдущие направления, что позволяет бороться с известной проблемой ``зигзага'' в методе FW \cite{Braun2022} аналогично методу сопряженных градиентов.

Алгоритм CFW хорошо описан в статье \cite{mitradjieva2012stiff}. В статье доказывается глобальная сходимость алгоритма CFW.  В статье также предлагается алгоритм Bi-conjugate Frank-Wolfe (BFW), в котором теперь берется выпуклая комбинация из двух предыдущих направлений и шага FW, а новое направление $d^k$ должно быть сопряженным к этим предыдущим направлениям.

В текущей статье предлагается алгоритм N-conjugate Frank-Wolfe (NFW), обобщающий CFW до числа сопряженных направлений равного $N > 2$. То есть используются 3 и более предыдущих направления для выбора следующего направления шага, сопряженного ко всем предыдущим N направлениям.



Другой известной модификацией FW является алгоритм Fukushima Frank-Wolfe (FFW), предложенный в статье \cite{FUKUSHIMA1984169}. Метод хорошо описан в статье \cite{arrache2008accelerating},  в которой также предлагаются некоторые модификации алгоритма. Данный алгоритм из того же класса, что и CFW, BFW, NFW, ведь использует для следующего направления шага предыдущие направления FW, а именно берет их среднее.

В текущей статье предлагается алгоритм-модификация Weighted Fukushima Frank-Wolfe (WFFW). 
Данный алгоритм использует вместо обычного среднего по направлениям~--- экспоненциальное сглаживание.

Перечислим ключевые результаты данной статьи:
\begin{enumerate}
    \item Дан обзор существующих модификаций алгоритма FW для задачи распределения потоков;
    \item Предложены новые модификации алгоритма FW: NFW и WFFW;
    \item Проведены численные эксперименты с использованием открытого датасета \cite{TNRT}. Установлено, что в большинстве случаев наиболее эффективным оказывается алгоритм NFW при $N=3$. 
\end{enumerate}

\paragraph{Постановка задачи}

В данной работе рассматривается модель Бекмана.
Задача равновесного распределения потоков сводится к минимизации следующей потенциальной функции: 
\begin{gather*}
    \Psi(f) = \sum_{e \in \mathcal{E}} \underbrace{\int_{0}^{f_e} \tau_e (z) d z}_{\sigma_e(f_e)} 
    \longrightarrow \min_{f = \Theta x, \; x \in X} \tag{Primal}.
\end{gather*}

В качестве функций затрат $\tau_e$, оценивающих время проезда по дорогам при заданных потоках, используются классические BPR-функции:
\begin{gather*}         
         \tau_e(f_e) = \bar{t}_e \left(1 + \rho \left( \frac{f_e}{\bar{f}_e}\right)^{\frac{1}{\mu}} \right) \tag{Beckman}. \\
         \hspace{1cm} \tau(f) \equiv \{\tau_e(f_e)\}_{e \in \mathcal{E}} - \textbf{Время проезда по дорогам с заданным потоком} \tag{times}.\\
    \hspace{1cm} \rho , \mu , \overline{t_e} , \overline{f_e} - \textbf{Константы в модели Бекмана}.
\end{gather*}

Задан граф дорог, а задача --- найти оптимальное распределение потоков: 
\begin{gather*}
    \mathcal{E} - \textbf{Ребра дорожного графа} \tag{roads}.\\
    \hspace{1cm} f \equiv \{f_e\}_{e \in \mathcal{E}} - \textbf{Потоки по дорогам}. \tag{flows}
\end{gather*}

Предполагается, что заданы пары источников и стоков, а также корреспонденции по ним:  
\begin{gather*}
    \hspace{1cm} OD  \equiv \{ (i,j) | i \in O , j \in D \} - \textbf{ Множество пар источник-сток} \tag{origin-destination}. \\
    \hspace{1cm} \{d_{(i,j)}\}_{(i,j) \in OD} - \textbf{ Корреспонденции или заданный поток машин из }i\textbf{ в }j \tag{correspondences}. \\
    \hspace{1cm} \mathcal{P} = \{\mathcal{P}_{i,j}\}_{(i,j)\in OD} - \textbf{Множество всех путей из O в D}. 
\end{gather*}

Корреспонденции, соответствующие каждой паре источник-сток, аддитивно разбиваются на потоки по путям, связывающим эти пары:
\begin{gather*}
    \hspace{1cm} X = \{ x \; | \; d_{(i,j)} = \sum\limits_{p\in \mathcal{P}_{i,j}} x_p \; , \; x_p \ge 0 \; \forall (i,j)\in OD  \}. 
\end{gather*}

Потоки по всем путям, проходящим по некоторой дороге, накладываются и порождают некоторую величину потока на этой дороге. Это можно записать компактно, введя матрицу принадлежности ребер путям:  

\begin{gather*}
    \hspace{1cm} \Theta \equiv \{\delta_{e,p}\}_{e \in \mathcal{E},p\in \mathcal{P}} - \textbf{Матрица принадлежности ребер} e \textbf{ путям } p . \\
    \hspace{1cm} x \equiv \{x_{p}\}_{p \in \mathcal{P}} - \textbf{Потоки по заданным путям } p\tag{flow on path}. \\
    \hspace{1cm} f = \Theta x - \textbf{Связь потоков по дорогам и по путям}. 
\end{gather*}

Сама задача подробно описана в работе \cite{KotlyarovaEtAl2021}.

\paragraph{Алгоритмы}

 В данном параграфе все модификации алгоритма Frank-Wolfe используют общий подход: поиск направления шага $d^k$ с использованием либо последних направлений FW $y_k^{FW}$, либо последних шагов $d^{k-m} \equiv s^{k-m} - f^{k-m}$.

\subparagraph{Frank-Wolfe} 

Алгоритм Frank-Wolfe можно разбить на два этапа: линейная аппроксимации целевой функции с последующей минимизацией линейного функционала и поиск направления шага с применением одномерной минимизации.

В разложении в ряд Тейлора целевой функции оставляется только первый линейный член, который потом используется в минимизации функции:
\begin{gather*}
    \Psi_k(y) \equiv \Psi(f^k) + \langle \nabla \Psi(f^k) , y - f^k \rangle \tag{linear approx}.
\end{gather*}

Решается линейная задача: 
\begin{gather*}
    \min\limits_{y = \Theta x, \; x \in X}\Psi_k(y) \iff \min\limits_{y = \Theta x, \; x \in X} \langle \nabla \Psi(f^k) , y \rangle \tag{linear problem (LP)}.
\end{gather*}

Далее делается шаг в направлении решения линейной задачи, а длина шага подбирается с помощью линейного поиска:  
\begin{gather*}
    s_k^{FW} \equiv \argmin\limits_{y = \Theta x, \; x \in X} \langle \nabla \Psi(f^k) , y \rangle = \argmin\limits_{y = \Theta x, \; x \in X} \langle \tau(f^k) , y \rangle \tag{solve LP}
    .\\
    d_k^{FW} = s_k^{FW} - f^k \tag{directional}. 
    \\
    f^{k+1} = f^k + \gamma_k d_k^{FW} , \;
    \gamma_k = \argmin_{\gamma \in [0,1]}(f^k + \gamma d^{FW}_k) \tag{\textbf{linesearch}}.
\end{gather*}

Стоит заметить, что наиболее вычислительно затратная операция в алгоритме FW для задачи распределения потоков~--- решение линейной подзадачи алгоритмами поиска кратчайших путей.
Сама же целевая функция представляется в виде суммы простых одномерных функций, и вычисление её и её производных значительно проще поиска кратчайших путей в графе.
Поэтому добавление в алгоритм FW процедур линейного поиска или описанного далее поиска сопряженных направлений не приводит к заметному увеличению сложности итерации.

В дальнейшем решение линейной задачи будет называться направлением FW и обозначаться $y_k^{FW}$ или $s_k^{FW}$. При этом данная линейная задача минимизации решается путем поиска кратчайших путей в графе дорог с текущими затратами $\tau(f^k)$, ведь $\langle \tau(f^k) ,y \rangle $ равно суммарным затратам при запуске потоков $y$ по дорогам. 

На данный момент известна сходимость алгоритма по прямой функции $\Psi(f)$ и по зазору FW из статьи \cite{pedregosa2020linearly}, где $L$~--- константа Липшица градиента, $\diam(C)$~--- диаметр множества $\{f| f = \Theta x , x \in X \}$. 

\begin{gather*}
   g_k \equiv \max\limits_{s = \Theta x , x \in X } \langle \nabla\Psi(f^k) , f^k - s \rangle  \tag{frank-wolfe gap}. \\
    g^*_k = \min\limits_{ m \in \overline{1..k}} g_m \tag{min FW gap}.\\
   g^*_k \le \frac{\min\{ 2h_0 , L \diam(C)^2\}}{\sqrt{k + 1}}\tag{FW gap convergence}. \\
   \Psi(f^k) - \Psi(f^*) \le \frac{2L\diam(C)^2}{k + 1} \tag{primal convergence}.
\end{gather*}

Эти результаты применимы и к FW в задаче поиска равновесного распределения потоков с моделью Бекмана. Ведь диаметр множества потоков по путям ограничен в силу того, что оно является прямым произведением многомерных симплексов. Липшицевость градиента обеспечивается на этом компакте.

\begin{algorithm}[H]
    \begin{algorithmic}[1]
    \STATE $t^0 \coloneqq \bar{t} \quad ,f^0 \coloneqq \argmin\limits_{s \in \{\Theta x : x \in X\}} \langle t^0, s \rangle$, $k \coloneqq 0$
    \REPEAT
        \STATE $s^k \coloneqq \argmin\limits_{s \in \{\Theta x : x \in X\}} \langle t^k, s \rangle \quad ,t_e^k \coloneqq \frac{\partial \Psi (f^k)}{\partial f_e} = \tau_e(f^k)$ 
        \STATE $\gamma_k \coloneqq \argmin\limits_{\gamma \in [0,1]} ((1 - \gamma) f^k + \gamma s^k) \quad or \quad \frac{2}{k + 1}$
        \STATE $f^{k + 1} \coloneqq (1 - \gamma_k) f^k + \gamma_k s^k$ 
        \STATE $k \coloneqq k + 1$
    \UNTIL{ k < max\_iter}
    \end{algorithmic}
    \caption{Frank--Wolfe algorithm}
    \label{alg::frank-wolfe}
\end{algorithm}

\subparagraph{Conjugate Frank-Wolfe}

Идея данного алгоритма заключается в выборе направления шага $s^k$, как выпуклой комбинации предыдущего направления $s^{k-1}$ и направления $s^{FW}_k$, соответствующего обычному FW шагу: 
\begin{gather*}
    s^k = \alpha_0 s_k^{FW} + \alpha_1 s^{k-1} \tag{convex combination}.
\end{gather*}

Выбор же коэффициентов в выпуклой комбинации находится из условия сопряженности нового шага $d^k$ к предыдущему $d^{k-1}$ по гессиану минимизируемого функционала $\nabla^2 \Psi(f^k)$ :
\begin{gather*}
    (d^k)^T\nabla^2 \Psi(f^k) d^{k-1} = 0 \tag{conjugate}.
\end{gather*}

Гессиан $\nabla^2 \Psi(f^k)$ находится аналитически в случае модели Бекмана и вычисляется за такое же время, как и градиент. 

В статье \cite{mitradjieva2012stiff} также предлагается алгоритм Bi-conjugate Frank-Wolfe (BFW), в котором теперь берется выпуклая комбинация из двух предыдущих направлений и шага FW, а новое направление $d^k$ должно быть сопряженным к этим предыдущим направлениям по гессиану в текущей точке $f^k$:
\begin{gather*}
    (d^k)^T\nabla^2 \Psi(f^k) d^{k-1} = 0. \\
        (d^k)^T\nabla^2 \Psi(f^k) d^{k-2} = 0 \tag{bi-conjugate}.
\end{gather*}
При этом делается предположение о сопряженности двух предыдущих направлений по новому гессиану, позволяющее получить коэффициенты выпуклой комбинации: 
\begin{gather*}
    (d^{k-1})^T\nabla^2 \Psi(f^k) d^{k-2} = 0 \tag{assumption}.
\end{gather*}

Вывод этих коэффициентов написан в статье \cite{mitradjieva2012stiff}.

\begin{algorithm}[H]
\begin{algorithmic}[2]
    \STATE $t^0 \coloneqq \bar{t} \quad ,f^0 \coloneqq \argmin\limits_{s \in \{\Theta x : x \in X\}} \langle t^0, s \rangle \quad , k \coloneqq 0 \quad, x^* \coloneqq f^0$
    \STATE FWM step for initialize $f^1$
    \REPEAT
        \STATE $s^k \coloneqq \argmin\limits_{s \in \{\Theta x : x \in X\}} \langle t^k, s \rangle \quad ,t_e^k \coloneqq \frac{\partial \Psi (f^k)}{\partial f_e} = \tau_e(f^k)$ 
        \STATE $\alpha \coloneqq \frac{(x^*-f^k)^TH^k(s^k-f^k)}{(x^*-f^k)^TH^k(s^k-f^{k-1})} $ \textbf{, где} $H^k = \nabla^2 \Psi(f^k)$ 
        \STATE Project $\alpha$ to [0,$\alpha_0$]
        \STATE $x^* \coloneqq \alpha x^* + (1 - \alpha) s^k $
        \STATE $\gamma_k \coloneqq \argmin\limits_{\gamma \in [0,1]} ((1 - \gamma) f^k + \gamma x^*)$
        \STATE $f^{k + 1} \coloneqq (1 - \gamma_k) f^k + \gamma_k x^*$ 
        \STATE $k \coloneqq k + 1$
    \UNTIL{ k < max\_iter}
    \end{algorithmic}
    \caption{Conjugate Frank--Wolfe algorithm}
    \label{alg::conjugate-frank-wolfe}
\end{algorithm}

\subparagraph{Fukushima Frank-Wolfe}

Автор статьи \cite{FUKUSHIMA1984169} предлагает брать среднее от последних $l$ направлений FW $\{ y_{k-m}^{FW}\}_{m=0}^l$:
\begin{gather*}
s^k = \frac{1}{l+1} \cbraces{s_k^{FW} + \sum^{l}_{i=1}s_{k-i}^{FW}}  \tag{mean directions}.
\end{gather*}

Также добавляется условие на оптимальность выбора этого среднего по сравнению с выбором только последнего направления FW.

\begin{algorithm}[H]
\begin{algorithmic}[3]
    \STATE $l \in \mathbb{N}$ 
    \STATE $t^0 \coloneqq \bar{t} \quad ,f^0 \coloneqq \argmin\limits_{s \in \{\Theta x : x \in X\}} \langle t^0, s \rangle$, $k \coloneqq 0$ 
    \REPEAT
        \STATE $s^k := \argmin\limits_{s \in \{\Theta x : x \in X\}} \langle t^k, s \rangle \quad , t_e^k \coloneqq \frac{\partial \Psi (f^k)}{\partial f_e} = \tau_e(f^k)$ 
        \STATE $q = \min( k , l ) - 1 $  
        \STATE $\nu_k = \sum\limits_{i = k - q}^{k} \mu_i s^i - f^k \quad ,\sum\limits_{i = k - q}^{k} \mu_i = 1 , \mu_i \ge 0 $
        \STATE $w_k = s^k - f^k$
        \STATE 

        \STATE $d^k =
            \begin{cases}
                \begin{aligned}
                    & \nu_k, \quad &&\text{если }   \frac{\nabla\Psi(f^k)^T\nu_k}{\|\nu_k\|} \le  \frac{\nabla\Psi(f^k)^Tw_k}{\|w_k\|} \\
                    &w_k, \quad 
                    &&\text{иначе } 
                \end{aligned}
            \end{cases}$
        \STATE $\gamma_k \coloneqq \argmin\limits_{\gamma \in [0,1]} (f^k + \gamma d^k)$
        \STATE $f^{k + 1} \coloneqq f^k + \gamma_k d^k$
        \STATE $k \coloneqq k + 1$
    \UNTIL{ k < max\_iter}
    \end{algorithmic}
    \caption{Fukushima Frank--Wolfe algorithm}
    \label{alg::fukushima-frank-wolfe}
\end{algorithm}

\subparagraph{Weighted Fukushima Frank-Wolfe}

Идея брать скользящее среднее по предыдущим направлениям в FFW привела к идее делать экспоненциальное сглаживание по этим направлениям. Поэтому предлагается Weighted FFW (WFFW), в котором производится экспоненциальное сглаживание по направлениям FW $s_k^{FW}$:

\begin{gather*}
    Q_{k} =  \beta Q_{k-1} + (1 - \beta) s_k^{FW} \tag{exponential smoothing}\\
    d_k = Q_{k} - f^k \tag{direction}.
\end{gather*}

\begin{algorithm}[H]
\begin{algorithmic}[4]
    \STATE $W \in [0,1] $ 
    \STATE $t^0 \coloneqq \bar{t} \quad ,f^0 \coloneqq \argmin\limits_{s \in \{\Theta x : x \in X\}} \langle t^0, s \rangle$, $k \coloneqq 0 \quad , Q^0 = f^0$ 
    \REPEAT
        \STATE $s^k := \argmin\limits_{s \in \{\Theta x : x \in X\}} \langle t^k, s \rangle \quad , t_e^k \coloneqq \frac{\partial \Psi (f^k)}{\partial f_e} = \tau_e(f^k)$ 
        \STATE  $Q^k = (1 - W)Q^{k-1} + Ws^k$
        \STATE $d^k = Q^k - f^k$

        \STATE $\gamma_k \coloneqq \argmin\limits_{\gamma \in [0,1]} (f^k + \gamma d^k)$
        \STATE $f^{k + 1} \coloneqq f^k + \gamma_k d^k$
        \STATE $k \coloneqq k + 1$
    \UNTIL{ k < max\_iter}
    \end{algorithmic}
    \caption{Weighted Fukushima Frank--Wolfe algorithm}
    \label{alg::weighted-fukushima-frank-wolfe}
\end{algorithm}

\newpage

\subparagraph{N-conjugate Frank-Wolfe}

Как и в случае Conjugate Frank-Wolfe берется выпуклая комбинация $s_k^{FW}$ и предыдущих направлений $\{s^{k-i}\}_{i=1}^{N_{curr}}$. Имеется степень свободы, связанная с выбором коэффициентов выпуклой комбинации. Эта степень свободы позволяет сделать следующее направление сопряженным к остальным, если выполняется предположение, что предыдущие направления сопряжены между собой по новому посчитанному гессиану:

\begin{gather*}
 (d^{k - m})^T \nabla^2 \Psi(f^k) d^{k - n} = 0 \quad \forall m \neq n \in \overline{1, N_{curr}} \tag{assumption}.
\end{gather*}
Оно техническое, ведь в дальнейшших выкладках упрощает расчет коэффициентов, зануляя слагаемые в сумме. Причем оно не лишено смысла: если на предыдущих итерациях направления сопряжены по предыдущим гессианам, а гессиан медленно меняется, то с некоторой точностью предыдущие направления сопряжены и по новому гессиану. 

Благодаря предположению, все направления, включая новое посчитанное, $\{ d^{k-j} \quad \forall j \in \overline{0, N_{curr}}\}$ будут сопряжены между собой попарно. 

Следующий шаг алгоритма~--- это выпуклая комбинация предыдущих направлений и нового посчитанного направления FW:
\begin{gather}
    s^k = \alpha_0 s_k^{FW} + \sum\limits_{m=1}^{N_{curr}} \alpha_m s^{k-m} \tag{convex combination}.
\end{gather}

Это можно переписать в терминах $d^k$:
\begin{gather}
    d^k = \alpha_0 s_k^{FW} + \sum\limits_{m=1}^{N_{curr}} \alpha_m s^{k-m} - f^k  \tag{direction}.
\end{gather}

Из определения $f^k = (1- \gamma_{k-1} )f^{k-1} + \gamma_{k-1} s^{k-1} $ и $d^{k-1} = s^{k-1} - f^{k-1} $ можно получить:
\begin{gather}
    \overline{d^{k-1}} \equiv (1 - \gamma_{k-1}) d^{k-1}. \\ 
    s^{k-1} - f^k = \overline{d^{k-1}}.
\end{gather}

Запишем телескопическую сумму и внедрим ее в сумму, задающую направление: 
\begin{gather}
    s^{k-m} - f^k = \sum\limits_{i = 0}^{ m -2}\cbraces{s^{k-m+i}  - s^{k - m + i + 1}} + s^{k-1} - f^k \tag{telescop sum}.\\
    d^k = \alpha_0 d_k^{FW} + \sum\limits_{m=1}^{N_{curr}} \alpha_m \cbraces{\sum\limits_{i = 0}^{ m - 2} \cbraces{s^{k-m+i}  - s^{k - m + i + 1}} + \overline{d^{k-1}}}.
\end{gather}

Перейдем от терминов $s^k$ к терминам $d^k$ и $\overline{d^k}$, воспользовавшись выражением (2), а после перепишем сумму в этих терминах и получим удобную запись направления шага NFW: 
\begin{gather}
    s^{k-m+i}  - s^{k - m + i + 1} = s^{k-m+i}  - s^{k - m + i + 1} - f^{k - m + i + 1} + f^{k - m + i + 1} = \overline{d^{k - m + i}} - d^{k - m + i + 1}. \\
    d^k = \alpha_0 d_k^{FW} + \sum\limits_{m=1}^{N_{curr}} \alpha_m \cbraces{\sum\limits_{i = 0}^{ m - 2} \cbraces{\overline{d^{k-m+i}}  - d^{k - m + i + 1}} + \overline{d^{k-1}}}. \\
    d^k = \alpha_0 d_k^{FW} + \sum\limits_{m=1}^{N_{curr}} \alpha_m \cbraces{\sum\limits_{i = 0}^{ m - 1} \overline{d^{k-m+i}}  - \sum\limits_{i = 0}^{ m - 2} d^{k - m + i + 1} } \tag{direction NFW}.  
\end{gather}

Теперь внедрим результаты выкладок выше в условие сопряженности ($H^k \equiv \nabla^2 \Psi(f^k)$):

\begin{gather}
    d^{k-m} H^k d^k \tag{conjugate} = 0 .\\
    d^{k-m} H^k \sbraces{\alpha_0 d_k^{FW} + \sum\limits_{m=1}^{N_{curr}} \alpha_m \cbraces{\sum\limits_{i = 0}^{ m - 1} \overline{d^{k-m+i}}  - \sum\limits_{i = 0}^{ m - 2} d^{k - m + i + 1} }} = 0.
\end{gather}

Добавим некоторые переобозначения и выпишем еще раз предположение:
\begin{gather}
    d^{k - m} H^k d^{k - n} = 0 \quad \forall m \neq n \in \overline{1, N_{curr}} \tag{assumption}. \\
    A_{k-m} \equiv d^{k-m} H^k d_k^{FW} \tag{A}.\\
    B_{k-m} \equiv d^{k-m} H^k d^{k-m} \tag{B}.
\end{gather}

Выражение (direction NFW) имеет  много нулевых слагаемых, исходя из предположения. Поэтому после упрощения выражения можно последовательно выражать каждый коэффициент $\alpha_j$, начиная с $\alpha_N$:

\begin{gather}
    d^{k-m} H^k \sbraces{\alpha_0 d^{FW}_k + \overline{d^{k-m}} \sum\limits_{i = m}^{N_{curr}}\alpha_i  - d^{k-m}\sum\limits_{i = m+1}^{N_{curr}}\alpha_i} = 0 \tag{simplified}. \\
    \alpha_0A_{k-m} + (1 - \gamma_{k-m})B_{k-m}\alpha_m + B_{k-m} \cbraces{\sum\limits_{i = m + 1}^{N_{curr}}\alpha_i (- \gamma_{k - m})} = 0  \tag{add A and B}.
\end{gather}

В итоге можно получить выражения для $\alpha_m$. Делается переобозначение $\alpha_m \equiv \beta_m \alpha_0$, после можно последовательно вычислять $\beta_m$, начиная с $m = N_{curr}$ , заканчивая $m = 1$, сохраняя в памяти и последовательно изменяя сумму $\sum\limits_{n = m + 1}^{N_{curr}} \beta_n$:

\begin{gather}
    \beta_{N_{curr}} = \frac{-A_{k-N_{curr}}}{B_{k-N_{curr}}(1-\gamma_{k-N_{curr}})}. \\
    \beta_m = \frac{-A_{k-m}}{B_{k-m}(1-\gamma_{k-m})} + \frac{\gamma_{k-m}}{1-\gamma_{k-m}} \sum\limits_{n = m + 1}^{N_{curr}} \beta_n \quad \forall m = N_{curr}-1,...,1  .\\
    \alpha_m \equiv \beta_m \alpha_0 \quad \forall m \in \overline{1,N_{curr}}. \\
    \alpha_0 = \frac{1}{1 + \sum\limits_{i=1}^{N_{curr}}\beta_i} \tag{convex combination}.
\end{gather}

Если $\gamma_k = 1$, то текущие потоки $f^k$ не учитываются в получении новых $f^{k+1}$. В случае FW $f^{k+1}$ будет находиться на вершине симплекса, т.е. не в относительной внутренности. Истинное решение $f^*$ зачастую лежит в относительной внутренности. Для избежания такой ситуации используется следующая эвристика: когда $\gamma_k \ge \gamma_{max}$ , то <<обнуляются>> накопленные $N_{curr}$ направлений и получение $N$ сопряженных направлений начинается с начала. 

\begin{algorithm}[H]
\begin{algorithmic}[5]
    \STATE $N \in \mathbb{N} \quad \gamma_{max} \in [0,1]$ 
    \STATE $t^0 \coloneqq \bar{t} \quad, f^0 \coloneqq \argmin\limits_{s \in \{\Theta x : x \in X\}} \langle t^0, s \rangle$, $k \coloneqq 0 \quad, N_{curr} = 1 $ 
    \REPEAT
        \STATE $s_k^{FW}:= \argmin\limits_{s \in \{\Theta x : x \in X\}} \langle t^k, s \rangle \quad , t_e^k \coloneqq \frac{\partial \Psi (f^k)}{\partial f_e} = \tau_e(f^k)$ 
        \STATE Вычислить $\{\alpha_i\}_{i=0}^{N_{curr}}$
        \STATE $s^k = \alpha_0 s_k^{FW} + \sum\limits_{i=1}^{N_{curr}} \alpha_i s^{k-i} $        
        \STATE $d^k = s^k - f^k$
        \STATE $\gamma_k \coloneqq \argmin\limits_{\gamma \in [0,1]} (f^k + \gamma d^k)$
        \STATE $f^{k + 1} \coloneqq f^k + \gamma_k d^k$
        \STATE $N_{curr} =
            \begin{cases}
                \begin{aligned}
                    & N_{curr} + 1, \quad &&\text{если } N_{curr} < N \\
                    & N_{curr}, \quad &&\text{если } N_{curr} = N \\
                    & 1, \quad &&\text{если } \gamma_k > \gamma_{max} \\
                \end{aligned}
            \end{cases}$
        \STATE $k \coloneqq k + 1$
    \UNTIL{ k < max\_iter}
\end{algorithmic}
    \caption{N-conjugate Frank--Wolfe algorithm}
    \label{alg::N-conjugate-frank-wolfe}
\end{algorithm}

\newpage
\paragraph{Эксперименты}

В качестве оценки эффективности алгоритмов взята величина относительного зазора, которая также использовуется в статьях \cite{mitradjieva2012stiff} и \cite{arrache2008accelerating}: 
\begin{gather*}
    RG_k \equiv \frac{UBD_k - BLB_k}{BLB_k} \tag{relative gap}.
\end{gather*}

Относительный зазор представляет собой долю разницы между верхней оценкой минимума прямой функции и её нижней оценкой относительно последней. Эта величина служит индикатором близости текущего значения прямой функции к её минимуму в данное время.

В качестве верхней оценки берется значение прямой функции в следующей точке:
\begin{gather*}
    UBD_k \equiv \Psi(f^{k+1}) \tag{upper bound}. 
\end{gather*}

Чтобы получить нижнюю оценку в относительном зазоре вводится следующая нижняя оценка:
\begin{gather*}    
    LB_k \equiv \min\limits_{y = \Theta x, \; x \in X}\{\Psi(f^k) + \langle \nabla \Psi(f^k) , y - f^k\rangle\} \tag{def lower bound}.
\end{gather*}

Ее можно представить в виде разности прямой функции и зазора Frank-Wolfe:
\begin{gather*}
    g^k \equiv \langle \nabla \Psi(f^k) , f^k - s^{FW}_k\rangle \tag{FW gap} . \\
    LB_k = \Psi(f^k) + \langle \nabla \Psi(f^k) , f^k - s^{FW}_k\rangle = \Psi(f^k) - g^k \tag{LB with FW gap}.
\end{gather*}

В качестве нижней оценки для относительного зазора берется наибольшая нижняя оценка среди ранее посчитанных:
\begin{gather*}
    BLB_k \equiv \max\limits_{m \in \overline{1,k}} LB_m \tag{best lower bound}.
\end{gather*}

С использованием относительного зазора возможно провести верхнюю оценку отношения между разностью значения прямой функции в следующей точке и её минимальным значением к самому минимуму функции. Это можно показать, воспользовавшись выпуклостью прямой функции:

\begin{gather*}
    \Psi(f^*) \ge \Psi(f^{k}) + \langle \nabla \Psi(f^{k}) , f^* - f^k \rangle \ge \min\limits_{y = \Theta x, \; x \in X}\{\Psi(f^k) + \langle \nabla \Psi(f^k) , y - f^k\rangle\} \equiv LB_k \tag{convex}. 
    \\
    \Psi(f^*) \ge BLB_k .\\
    \frac{\Psi(f^{k+1}) - \Psi(f^*)}{\Psi(f^*)} \le \frac{\Psi(f^{k+1}) - BLB_k}{BLB_k} \equiv RG_k .
\end{gather*}

Также можно заметить, что данная оценка более точная, чем относительный зазор Frank-Wolfe: 
\begin{gather*}
    RG_k \le \frac{\Psi(f^k) - LB_k}{BLB} = \frac{g^k}{BLB}.
\end{gather*}

\subparagraph{Сравнение 3-conjugate FW с Bi-conjugate FW}
В данном разделе уделим внимание сравнению Bi-conjugate Frank-Wolfe с его обобщением N-conjugate Frank-Wolfe, так как в статье \cite{mitradjieva2012stiff} показано, что BFW эффективней чем CFW и FW. Но для того, чтобы показать превосходство алгоритмов NFW и BFW, в эксперименты также добавлены алгоритмы: FW, WFFW, CFW.

\begin{figure}[!ht]
\centering
\includegraphics[width=1.0\textwidth]{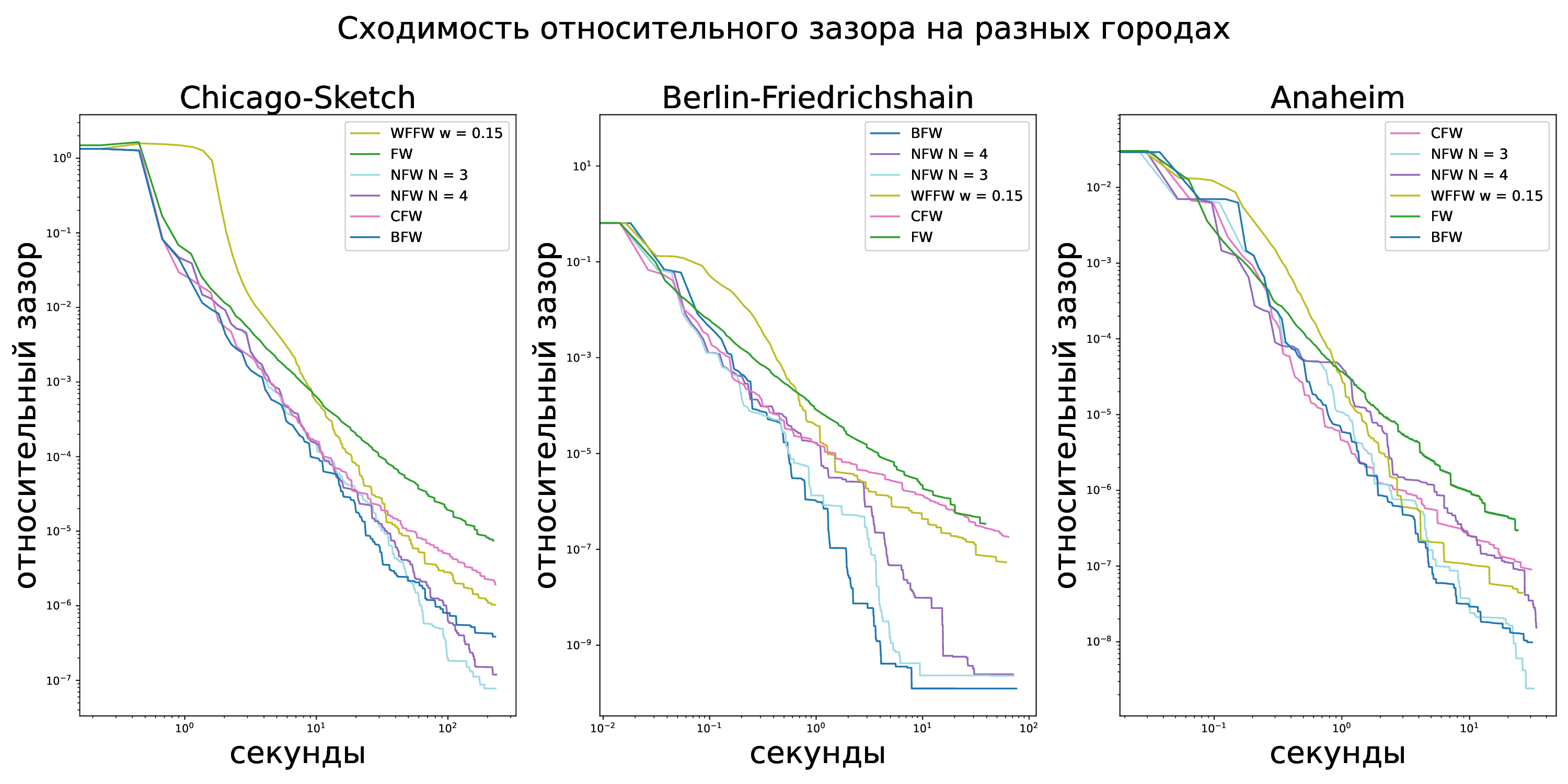}
\smallskip
\caption{Сходимость относительного зазора для алгоритмов NFW, CFW, BFW, WfFW, FW по времени (в секундах) на датасетах: Chicago-Sketch, Berlin-Friedrichshain, Anaheim. }
\label{fig:3FW_1}
\end{figure}

\begin{figure}[!ht]
\centering
\includegraphics[width=1.0\textwidth]{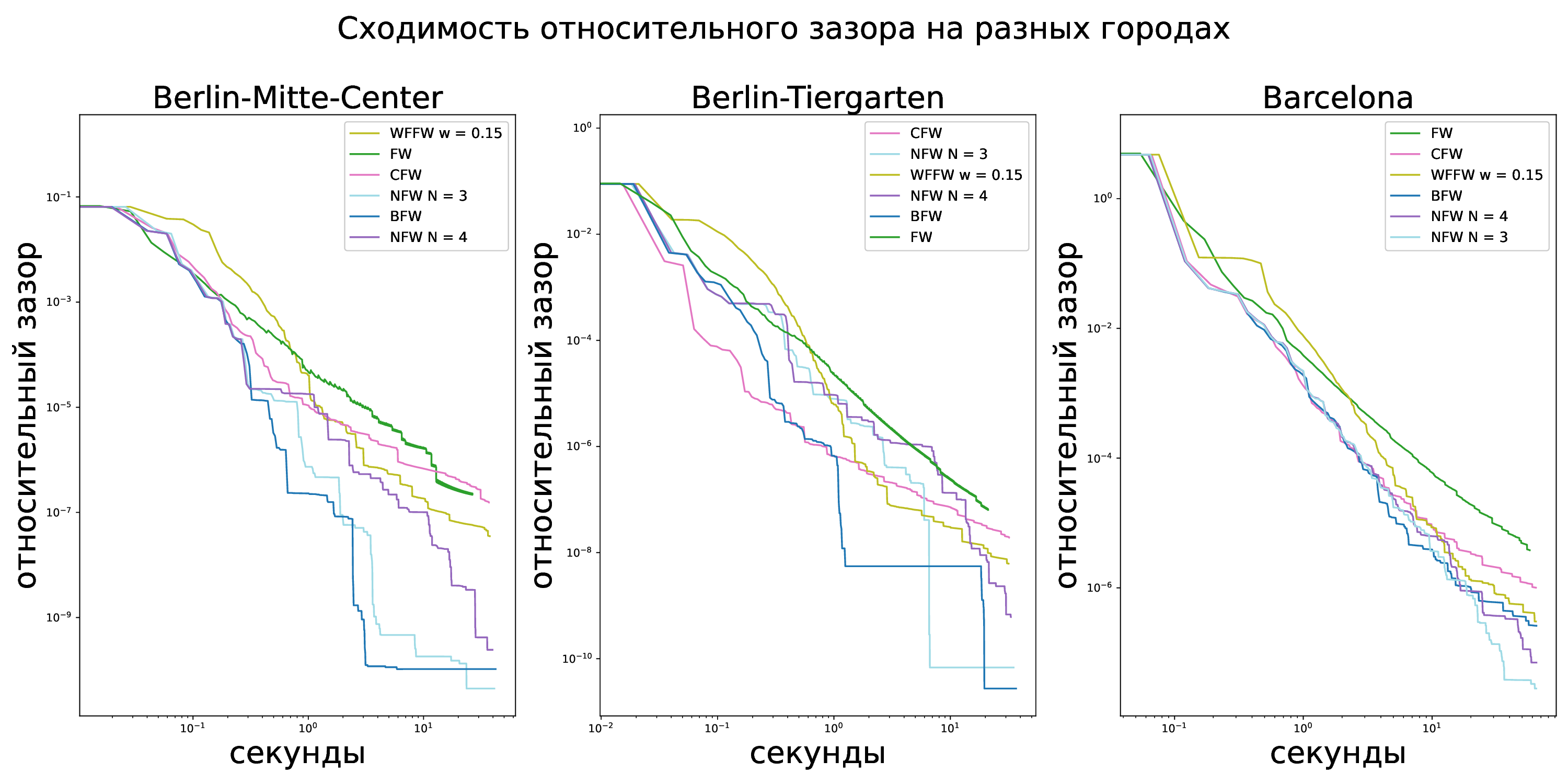}
\smallskip
\caption{Сходимость относительного зазора для алгоритмов NFW, CFW, BFW, WFFW, FW по времени (в секундах) на датасетах: Berlin-Mitte-Center, Berlin-Tiergarten, Barcelona. }
\label{fig:3FW_2}
\end{figure}

\begin{figure}[!ht]
\centering
\includegraphics[width=1.0\textwidth]{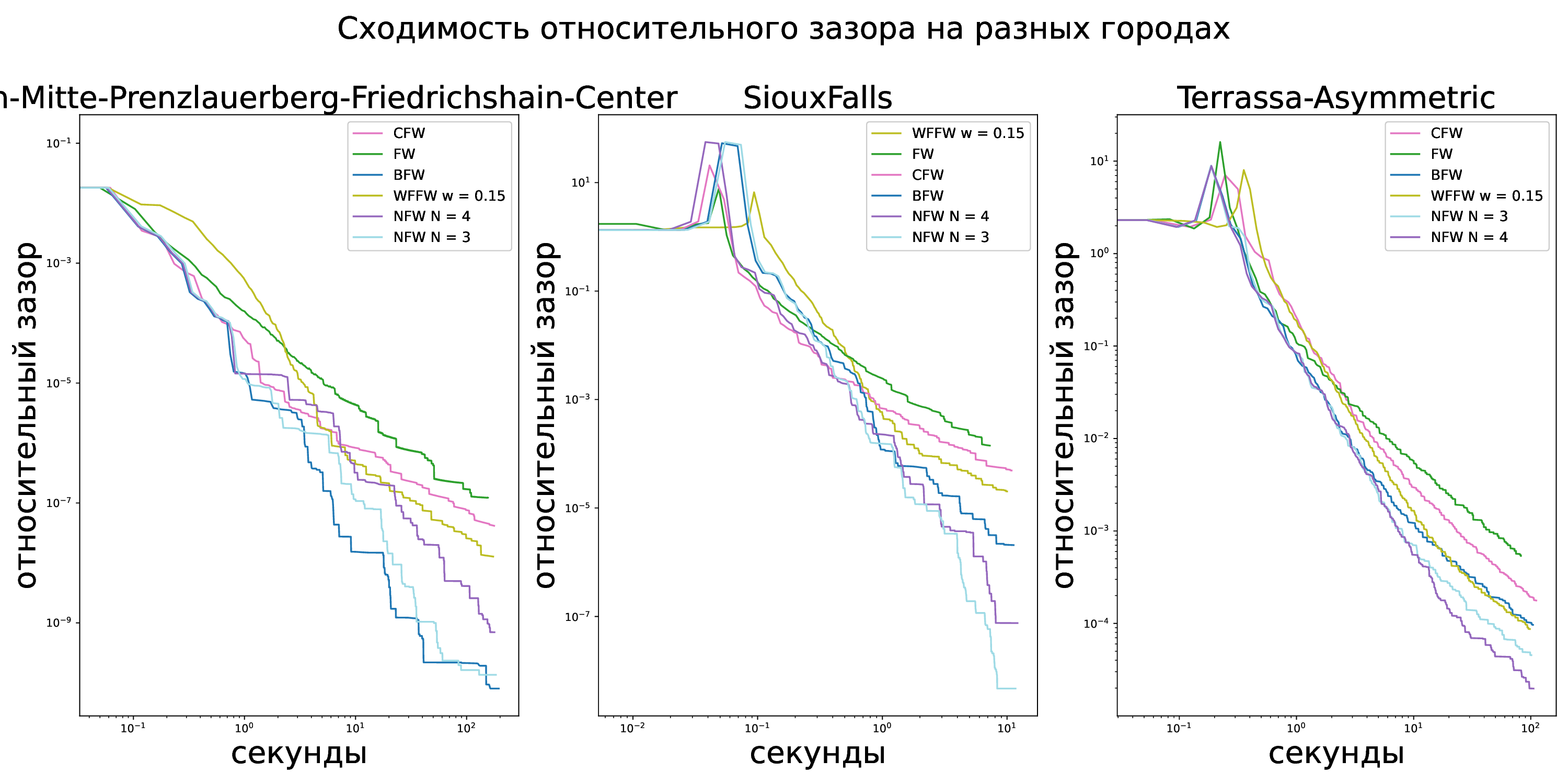}
\smallskip
\caption{Сходимость относительного зазора для алгоритмов NFW, CFW, BFW, WFFW, FW по времени (в секундах) на датасетах: Berlin-Mitte-Prenzlauerberg-Friedrichshain-Center, SiouxFalls, Terrassa-Asymmetric. }
\label{fig:3FW_3}
\end{figure}

Эксперименты были проведены на датасетах: Anaheim, SiouxFalls, Berlin-Tiergarten, Terrassa-Asymmetric, Chicago-Sketch, Berlin-Mitte-Center, Berlin-Friedrichshain, Barcelona, Berlin-Mitte-Prenzlauerberg-Friedrichshain-Center. Результаты изображены на Рис.~\ref{fig:3FW_1}, \ref{fig:3FW_2}, \ref{fig:3FW_3}. 

В итоге, в 6 из 9 случаев NFW с N = 3, побеждает BFW, а в остальных случаях проигрывает. Также можно заметить, что методы CFW, FW значительно уступают по эффективности остальным рассмотренным методам.

\subparagraph{Сравнение FFW с Weighted FFW}

\begin{figure}[!ht]
\centering
\includegraphics[width=1.0\textwidth]{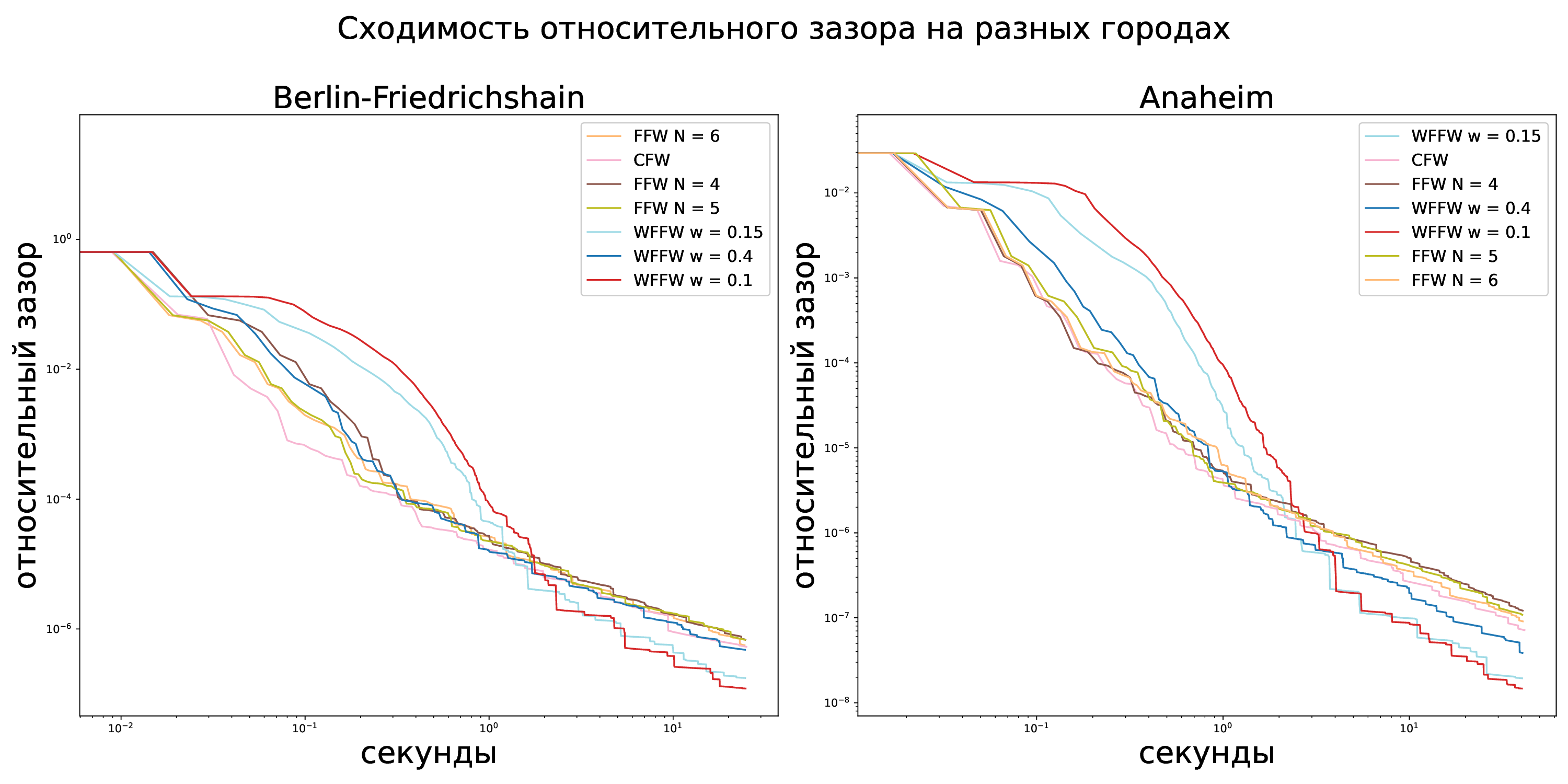}
\smallskip
\caption{Сходимость относительного зазора для алгоритмов FFW, WFFW, CFW по времени (в секундах) на датасетах: Anaheim, Berlin-Friedrichshain }
\label{fig:FFW_1}
\end{figure}

\begin{figure}[!ht]
\centering
\includegraphics[width=1.0\textwidth]{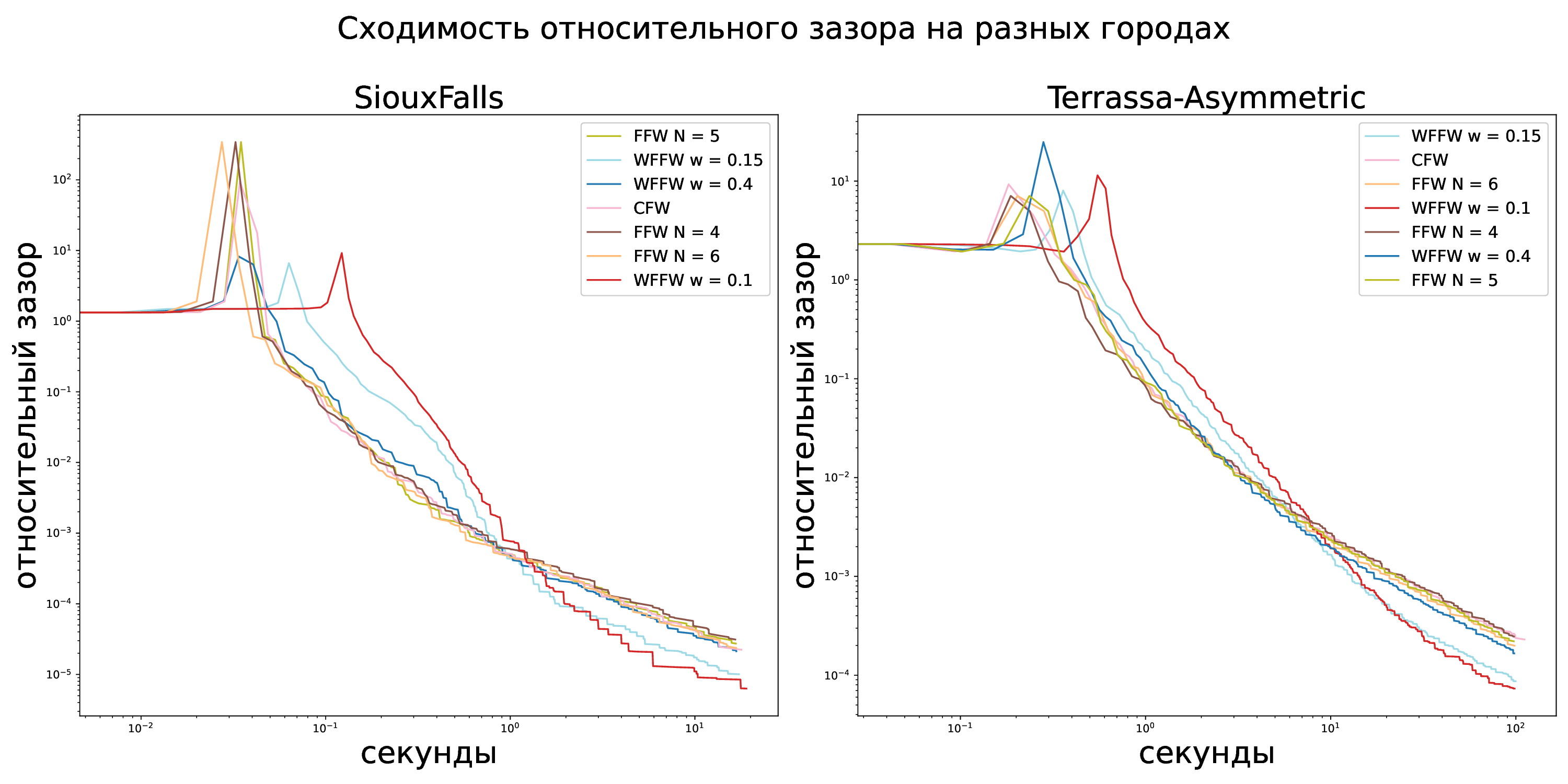}
\smallskip
\caption{Сходимость относительного зазора для алгоритмов FFW, WFFW, CFW по времени (в секундах) на датасетах: SiouxFalls, Terrassa-Asymmetric }
\label{fig:FFW_2}
\end{figure}

В данном разделе проведен сравнительный анализ не только между алгоритмом FFW и его модификацией, но также включен в рассмотрение CFW, поскольку CFW выступает в качестве базовой модификации алгоритма. Результаты приведены на Рис.~\ref{fig:FFW_1}, \ref{fig:FFW_2}.

Отмечается, что CFW проявляет сходство в процессе сходимости с FFW, в то время как WFFW демонстрирует значительное превосходство над обоими. Также можно заметить характерное поведение WFFW: в начальный момент времени он отстает от конкурирующих методов, набирая скорость, после чего обгоняет их в дальнейшем.

\subparagraph{Сравнение NFW для разных N}

\begin{figure}[!ht]
\centering
\includegraphics[width=1.0\textwidth]{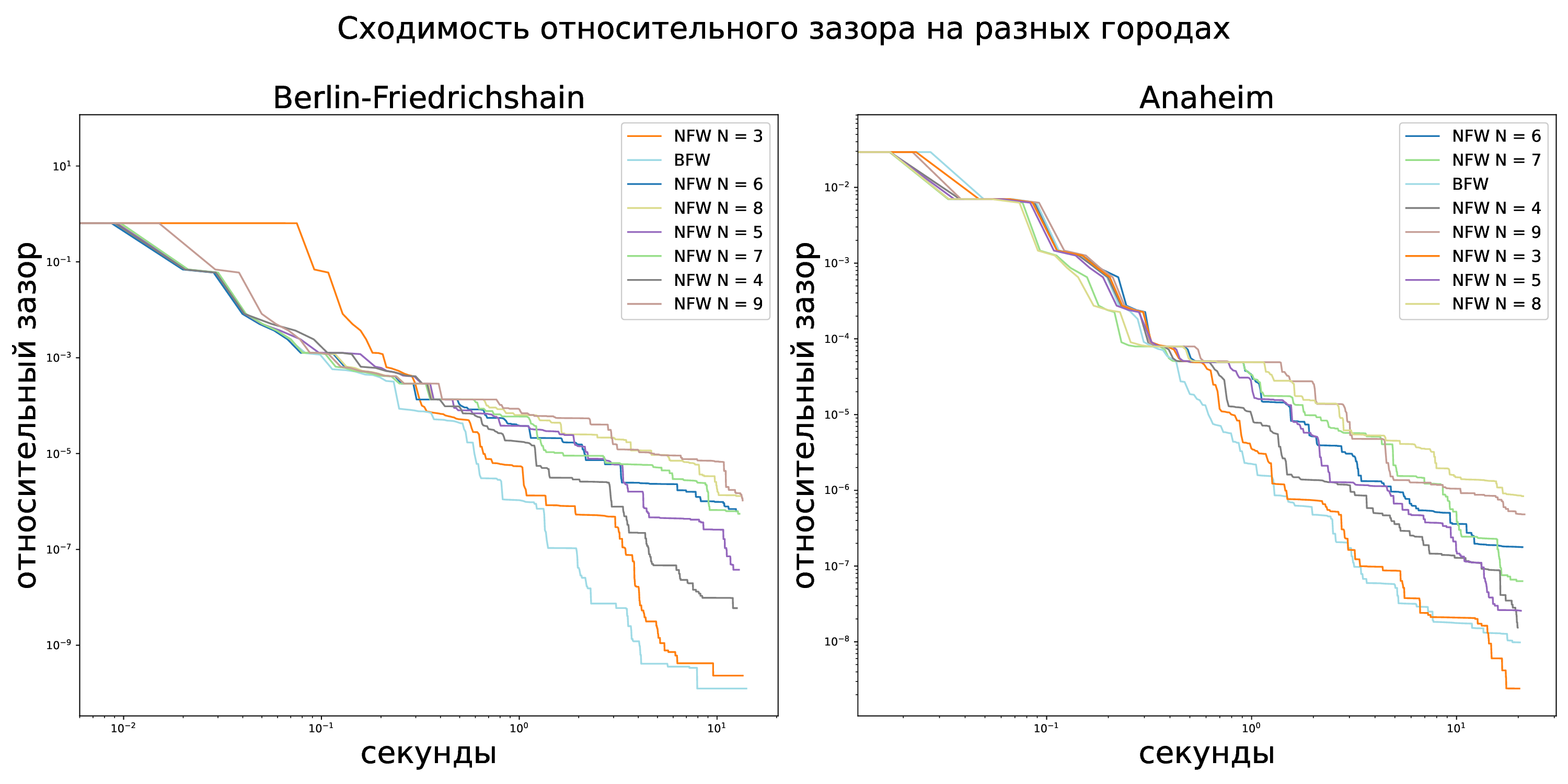}
\smallskip
\caption{Сходимость относительного зазора для алгоритмов NFW, BFW (в секундах) на датасетах: Anaheim, Berlin-Friedrichshain }
\label{fig:NFW_1}
\end{figure}

\begin{figure}[!ht]
\centering
\includegraphics[width=1.0\textwidth]{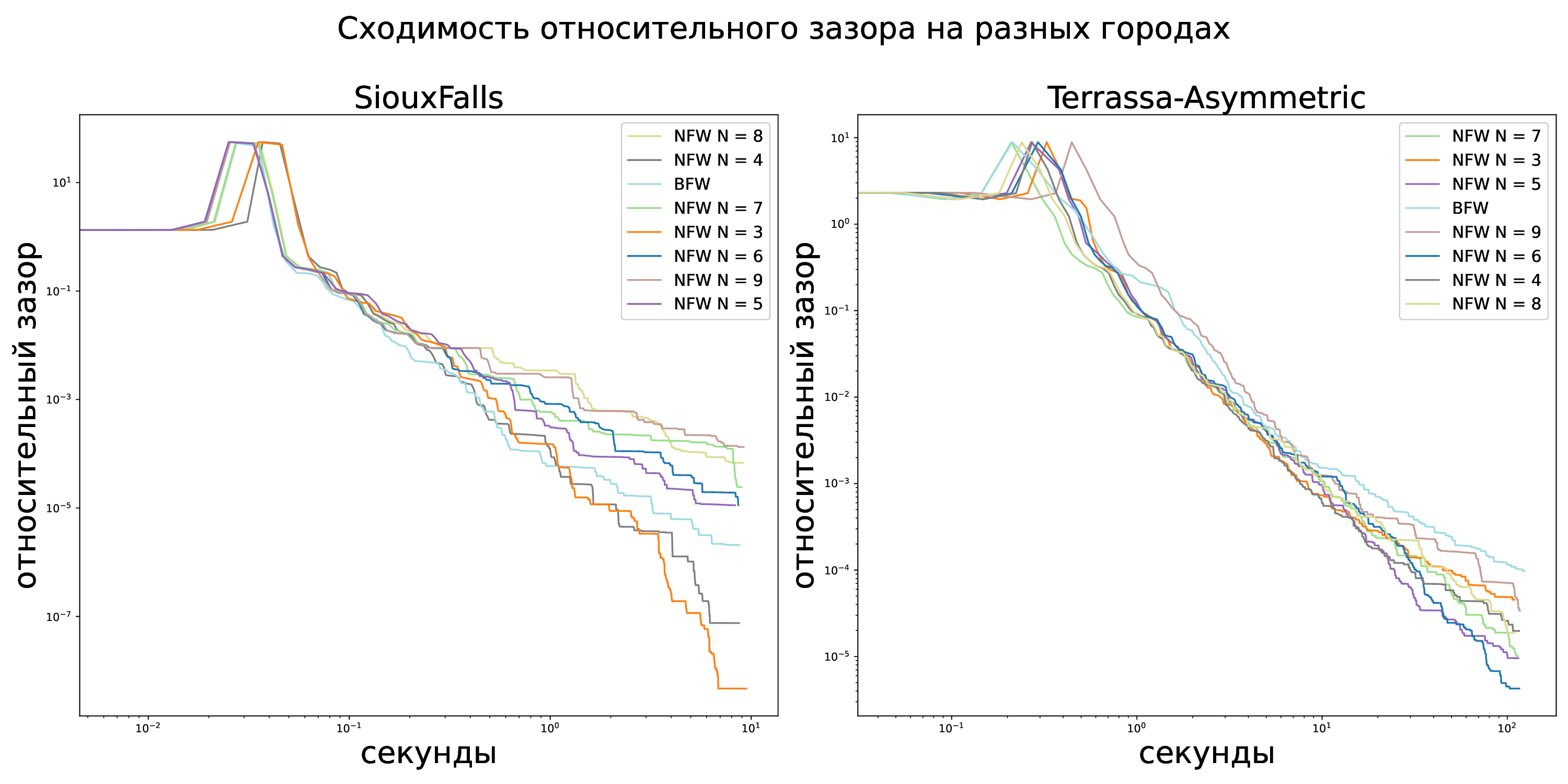}
\smallskip
\caption{Сходимость относительного зазора для алгоритмов NFW, BFW (в секундах) на датасетах: SiouxFalls, Terrassa-Asymmetric}
\label{fig:NFW_2}
\end{figure}

Также становится интересно проверить: как на эффективность алгоритмов NFW влияет гиперпараметр N. Результаты экспериментов приведены на Рис.~\ref{fig:NFW_1} , \ref{fig:NFW_2}.

На датасете SiouxFalls лучшей моделью оказалась модель с N = 3. На датасете Berlin-Friedrichshain, в свою очередь, лидировал Bi-conjugate Frank-Wolfe. На датасете Terrassa-Asymmetric лучшей моделью оказалась модель с N = 6. Интересным является неожиданное обстоятельство, заключающееся в том, что на данном датасете практически все модификации при различных значениях N продемонстрировали превосходство над моделью с N = 3.
Более того, Bi-conjugate FW показал себя менее эффективным, отстающим от всех рассмотренных алгоритмов.

Из результатов экспериментов можно сделать вывод, что при использовании гиперпараметра модели N = 3 не всегда достигаются лучшие результаты по сравнению с другими значениями N. Однако, в большинстве датасетов, алгоритм все же превосходит другие модификации. Интересно отметить, что иногда модели с N > 3 демонстрируют лучшие результаты, в то время как в некоторых случаях побеждает Bi-conjugate Frank-Wolfe.

\subparagraph{Сравнение WFFW для разных параметров}

\begin{figure}[!ht]
\centering
\includegraphics[width=1.0\textwidth]{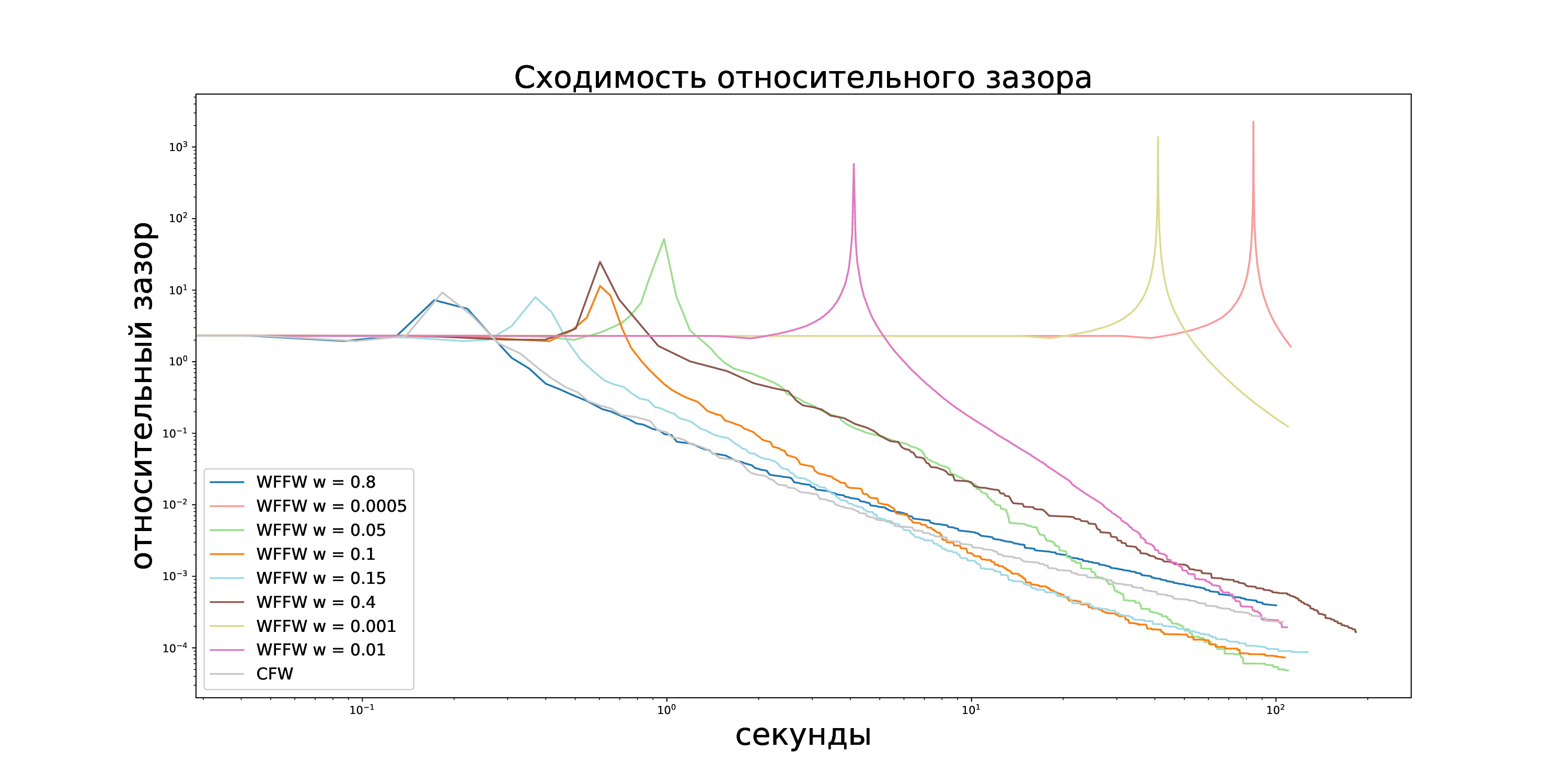}
\smallskip
\caption{Сходимость относительного зазора для алгоритмов WFFW с различными весами затухания (в секундах) на датасете: Terrassa-Asymmetric }
\label{fig:WFFW_1}
\end{figure}

На Рис.~\ref{fig:WFFW_1} изображены графики сходимостей WFFW при различных гиперпараметрах.

В результате анализа выявлен оптимальный выбор коэффициента затухания для заданного временного интервала. Выбор небольшого коэффициента подразумевает, что модель в значительной степени ориентируется на предыдущие направления и производит небольшие коррекции с их учетом, используя направление FW. Установлено, что уменьшение значения коэффициента приводит к увеличению времени, в течение которого относительный зазор незначительно измененяется со временем. После накопления скорости алгоритм в первую очередь движется в направлении, противоположном оптимальному направлению. Позже направление начинает настраиваться на направления FW, становясь оптимальным. Тем не менее, несмотря на данное простаивание, модель начинает догонять модели с более высокими значениями коэффициента. Однако, при использовании очень малых значений коэффициента, наблюдается значительное увеличение времени простаивания.

\paragraph{Вывод}

В настоящей статье были представлены разнообразные модификации алгоритма FW, подвергнуты сравнительному анализу на различных датасетах, представляющих реальные городские условия. В свете полученных результатов можно сделать вывод о том, что среди всех представленных модификаций наилучшие характеристики проявляют модели, относящиеся к классу N-conjugate Frank-Wolfe. Замечено, что часто оптимальным выбором является модель с N = 3, что соответствует использованию трех последних направлений в процессе оптимизации. Также было показано абсолютное превосходство алгоритма Weighted FFW перед алгоритмом FFW.

\end{document}